\input amstex
\documentstyle{amsppt}
%
%
\nopagenumbers
\def\Ker{\operatorname{Ker}}
\def\Cl{\operatorname{Cl}}
\hyphenation{ma-ni-fold}
\def\negskp{\hskip -2pt}
\def\compos{\,\raise 1pt\hbox{$\sssize\circ$} \,}
\pagewidth{360pt}
\pageheight{606pt}
\leftheadtext{Ruslan A. Sharipov}
\rightheadtext{Second problem of globalization \dots}
\topmatter
\title Second problem of globalization\\
{\lowercase{in the theory of dynamical systems
\vadjust{\vskip -1.0ex}admitting\\ the normal
shift of hypersurfaces.\vadjust{\vskip -2ex}}}
\endtitle
\author
R.~A.~Sharipov
\endauthor
\abstract
Problem of global integration of geometric structures arising
in the theory of dynamical systems admitting the normal shift
is considered. In the case when such integration is possible
the problem of globalization for shift maps is studied.
\endabstract
\address Rabochaya street 5, 450003, Ufa, Russia
\endaddress
\email \vtop to 20pt{\hsize=280pt\noindent
R\_\hskip 1pt Sharipov\@ic.bashedu.ru\newline
ruslan-sharipov\@usa.net\vss}
\endemail
\urladdr
http:/\negskp/www.geocities.com/CapeCanaveral/Lab/5341
\endurladdr
\endtopmatter
\loadbold
\document
\head
1. Introduction.
\endhead
     Let $M$ be Riemannian manifold of the dimension $n$. Newtonian
dynamical system in $M$ in local coordinates is determined by a system
of $n$ ODE's
$$
\ddot x^k+\sum^n_{i=1}\sum^n_{j=1}\Gamma^k_{ij}\,\dot x^i\,\dot x^j
=F^k(x^1,\ldots,x^n,\dot x^1,\ldots,\dot x^n),\hskip -2em
\tag1.1
$$
where $k=1,\,\ldots,\,n$. Here $\Gamma^k_{ij}=\Gamma^k_{ij}(x^1,
\ldots,x^n)$ are components of metric connection and $F^k$ are
components of force vector $\bold F$. They determine force field
of dynamical system \thetag{1.1}.\par
    The theory of Newtonian dynamical systems admitting the normal
shift of hypersurfaces was constructed in papers \cite{1--16}; on the
base of these papers two theses \cite{17} and \cite{18} were prepared.
We shall consider some details of this theory a little bit later. Now
we note only that this theory describes special class of force fields,
which in the case of higher dimensions $n\geqslant 3$ locally (in some
neighborhood of any point $p\in M$) can be given by explicit formula
$$
F_k=\frac{h(W)\,N_k}{W_v}-v\sum^n_{i=1}\frac{\nabla_iW}
{W_v}\,\bigl(2\,N^i\,N_k-\delta^i_k\bigr),\hskip -2em
\tag1.2
$$
where $W=W(x^1,\ldots,x^n,v)$ is some function of $(n+1)$ variables
and $h(w)$ is some function of one variable. The variable $v$ in
\thetag{1.2} denotes the modulus of velocity vector: $v=|\bold v|$.
While by $N^i$ and $N_k$ we denote contravariant and covariant
components of unitary vector $\bold N$ directed along velocity
vector:
$$
\xalignat 2
&N^i=\frac{v^i}{|\bold v|},&&N_k=\frac{v_k}{|\bold v|}.
\endxalignat
$$
Function $W(x^1,\ldots,x^n,v)$ in formula \thetag{1.2} should satisfy
the condition
$$
W_v=\frac{\partial W}{\partial v}\neq 0.\hskip -2em
\tag1.3
$$
This is quite natural, since partial derivative $W_v$ is in denominators
of two fractions in formula \thetag{1.2}.\par
     Suppose that Riemannian manifold $M$ is equipped with some Newtonian
dynamical system admitting the normal shift of hypersurfaces. Functions
$W$ and $h$ determine force field $\bold F$ of such system locally in
a neighborhood of some point $p\in M$. In the neighborhood of another
point $\tilde p\in M$ force field $\bold F$ in general case is determined
by another pair of functions $\tilde W$ and $\tilde h$. In the region of
overlapping of two neighborhoods (if they do really overlap) the force
field $\bold F$ can be determined by each of these two pairs of functions.
This gives an idea that force field $\bold F$ is related to some global
geometric structures on $M$, which are locally represented by pairs of
functions $(h,W)$. The problem of revealing such structures was called
{\bf the first problem of globalization}. It is solved by the following
theorem from \cite{19}.
\proclaim{Theorem 1.1} Defining Newtonian dynamical system admitting
the normal shift in Riemannian manifold $M$ is equivalent to defining
closed global section $\sigma$ for projectivized cotangent bundle
$\Cal P^*\!\Cal M$, where $\Cal M=M\times \Bbb R^{\sssize +}$,
satisfying the condition $\Ker\sigma\nparallel\Bbb R^{\sssize +}$,
and normalizing global section $s$ for one-dimensional factor-bundle
$\varOmega\Cal M=\Cal T\Cal M/U$, where $U=\Ker\sigma$.
\endproclaim
     Any section $\sigma$ of the bundle $\Cal P^*\!\Cal M$ in the
neighborhood of each point $q=(p,v)$ of $\Cal M=M\times\Bbb R^{\sssize
+}$ is determined by some differential 1-form $\boldsymbol\omega$,
which is unique up to scalar factor $\boldsymbol\omega\to\varphi\cdot
\boldsymbol\omega$. Closedness of $\sigma$ means that the form $\boldsymbol
\omega$ can be chosen closed. Each closed 1-form is locally exact,
it is a differential of some function:
$$
\boldsymbol\omega=dW\text{, \ where \ }W=W(x^1,\ldots,x^n,v).
$$
The condition $\Ker\sigma\nparallel\Bbb R^{\sssize +}$ means that
$(n+1)$-th component of the form $\boldsymbol\omega$ in local coordinates
$x^1,\,\ldots,\,x^n,\,v$ is nonzero:
$$
W_v=\omega_{n+1}\neq 0.\hskip -2em
\tag1.4
$$
In other words, $\Ker\sigma\nparallel\Bbb R^{\sssize +}$ is simply an
invariant (non-coordinate) form of the condition \thetag{1.3}. When
condition \thetag{1.4} is fulfilled, we can consider the quotients
$$
b_i=-\frac{\nabla_iW}{W_v}=-\frac{\omega_i}{\omega_{n+1}}.\hskip -2em
\tag1.5
$$
The quantities $b_1,\,\ldots,\,b_n$ do not change if we replace 
$\boldsymbol\omega$ by $\varphi\cdot\boldsymbol\omega$, they are
local coordinates in fibers of projectivized cotangent bundle
$\Cal P^*\!\Cal M$. If a section $\sigma$ of this bundle is given,
we have $n$ functions $b_i(x^1,\ldots,x^n,v)$, where $i=1,\,\ldots,\,n$.
The condition of closedness for $\sigma$ is written in form of the
following relationships:
$$
\left(\frac{\partial}{\partial x^j}+b_j\,\frac{\partial}
{\partial v}\right)b_i=\left(\frac{\partial}{\partial x^i}
+ b_i\,\frac{\partial}{\partial v}\right)b_j.\hskip -2em
\tag1.6
$$\par
    The kernel $\Ker\sigma=\Ker\boldsymbol\omega$ determines
$n$-dimensional distribution $U$ in the manifold $\Cal M=M\times
\Bbb R^{\sssize +}$, whose dimension is $n+1$. It also determines
$1$-dimensional vector-bundle obtained by factorization of
cotangent bundle $\Cal T\Cal M$ with respect to $U$. Due to the
condition $\Ker\sigma\nparallel\Bbb R^{\sssize +}$ vector field
$\partial/\partial v$ is transversal to $U$. Therefore each
section $s$ of factor-bundle $\varOmega\Cal M=\Cal T\Cal M/U$
in local coordinates can be determined by vector field of the
form
$$
a(x^1,\ldots,x^n,v)\cdot\frac{\partial}{\partial v},\hskip -2em
\tag1.7
$$
or by function $a(x^1,\ldots,x^n,v)$, which arises as a coefficient
in formula \thetag{1.7}. The concept of normalizing section is introduced
by the following two definitions from paper \cite{19}.
\definition{Definition 1.1} Vector field $\bold X$ is called {\bf
normalizing field} for smooth distribution $U$ if for any vector
field $\bold Y$ belonging to $U$ the commutator $[\bold X,\,\bold Y]$
is also in $U$.
\enddefinition
\definition{Definition 1.2} Section $s$ of factor-bundle $\varOmega\Cal M=
\Cal T\Cal M/U$ is called {\bf normalizing section} if in the neighborhood
of each point $q\in\Cal M$ it is represented by some normalizing vector
field for the distribution $U$.
\enddefinition
    The fact that section $s$ of factor-bundle $\Cal T\Cal M/\Ker\sigma$
is normalizing is expressed by the following equations for the function
$a$ in \thetag{1.7}:
$$
\left(\frac{\partial}{\partial x^i}+b_i\,\frac{\partial}{\partial v}
\right)a=\frac{\partial b_i}{\partial v}\,a.
\tag1.8
$$
Note that the concept of normalizing section $s$ of the bundle $\Cal T
\Cal M/U$ is correctly determined only for involutive distribution $U$.
In this case $s$ is a coset of vector $\bold X$ respective to subspace
$U$, i\.~e\. $s=\Cl_U(\bold X)$. The choice of vector field representing
such coset doesn't matter, since if $\bold X$ is normalizing vector field
for $U$ and $\bold Y\in U$, then the sum $\bold X+\bold Y$ is also
normalizing vector field for $U$. In our case $U=\Ker\sigma$ is involutive.
This follows from closedness of $\sigma$.
\head
2. Integration of geometric structures.
\endhead
     Theorem~1.1 determines global geometric structures related to
force fields of Newtonian dynamical systems admitting the normal shift,
thus solving first problem of globalization. As for calculation of
components of force vector, it yields formula
$$
F_k=a\,N_k+v\sum^n_{i=1}b_i\,\bigl(2\,N^i\,N_k-\delta^i_k\bigr),
\hskip -2em
\tag2.1
$$
where $a$ and $b_1,\,\ldots,\,b_n$ should satisfy the equations
\thetag{1.6} and \thetag{1.8}. It's clear that formula \thetag{2.1}
is much less effective than formula \thetag{1.2}. The passage from
\thetag{2.1} to \thetag{1.2} consists in integrating the equations
\thetag{1.6} and \thetag{1.8}. These equations are compatible and
locally integrable, this is shown in Chapter~\uppercase
\expandafter{\romannumeral 7} of thesis \cite{17} (see also paper
\cite{19}). Here we are interested in those cases, when they are
globally integrable.\par
     First step in global integration of the structures $\sigma$
and $s$ from theorem~1.1 consist in exploiting the closedness of
the section $\sigma$. On the manifold $\Cal M=M\times\Bbb R^{\sssize
+}$ (or possibly on universal cover for $\Cal M$) one should
find global closed 1-form $\boldsymbol\omega$ that should satisfy
the condition $\Ker\boldsymbol\omega=\Ker\sigma$. This form would
determine the quantities $b_1,\,\ldots,\,b_n$ according to the
formula \thetag{1.5}. If such formula is found, we say that
{\bf first level of global integration} of structures $\sigma$
and $s$ is reached.\par
    Assuming that first level of global integration of structures
$\sigma$ and $s$ is already reached, let's integrate the 1-form
$\boldsymbol\omega$ just found along the path binding some fixed
initial point $q_0$ with ending point $q$: 
$$
W(q)=\int^{\,q}_{q_0}\boldsymbol\omega.\hskip -2em
\tag2.2
$$
Passing from $\Cal M=M\times\Bbb R^{\sssize +}$ to universal cover
for $\Cal M$, if necessary, we warranty that formula yields single-valued
function $W$ on such cover. Then $dW=\boldsymbol\omega$.\par
    Let $\widehat M$ be universal cover for $M$. Then universal cover
$\widehat{\Cal M}$ for the manifold $\Cal M=M\times\Bbb R^{\sssize +}$
can be identified with $\widehat M\times\Bbb R^{\sssize +}$. The structure
of Cartesian product in $\widehat{\Cal M}=\widehat M\times\Bbb
R^{\sssize +}$ provides vector field $\bold V=\partial/\partial v$
directed along linear rulings in this manifold. Applying this field
to the function \thetag{2.2} we get the function $W_v=\bold VW$.
In local coordinates $x^1,\,\ldots,\,x^n,\,v$ this function coincides
with partial derivative:
$$
W_v=\frac{\partial W}{\partial v}=\omega_{n+1}\neq 0.
$$
Let's define the function $\widetilde W=a\cdot W_v$ and let's calculate
its differential:
$$
d\widetilde W=\sum^n_{i=1}\frac{\partial\widetilde W}{\partial x^i}
\cdot dx^i+\frac{\partial\widetilde W}{\partial v}\cdot dv.\hskip -2em
\tag2.3
$$
For first $n$ components in 1-form \thetag{2.3} we have
$$
\gather
\frac{\partial\widetilde W}{\partial x^i}=
\frac{\partial a}{\partial x^i}\,\omega_{n+1}+a\,\frac{\partial
\omega_{n+1}}{\partial x^i}=\left(\frac{\partial b_i}{\partial v}
\,a-b_i\,\frac{\partial a}{\partial v}\right)\omega_{n+1}+
a\,\frac{\partial\omega_{n+1}}{\partial x^i}=\\
=\left(-\frac{\partial\omega_i}{\partial v}\,a-b_i\,a\,
\frac{\partial\omega_{n+1}}{\partial v}-b_i\,\frac{\partial a}
{\partial v}\,\omega_{n+1}\right)+a\,\frac{\partial\omega_{n+1}}
{\partial x^i}=\frac{\omega_i}{\omega_{n+1}}\,\frac{\partial
\widetilde W}{\partial v}.
\endgather
$$
In these calculations we used closedness of the form $\boldsymbol\omega$
and the relationships \thetag{1.5} and \thetag{1.8}. The result of
calculations can be formulated as follows: the ratio of $i$-th and
$(n+1)$-th components of 1-form \thetag{2.3} is equal to the ratio of
$\omega_i$ and $\omega_{n+1}$. This means that forms $d\widetilde W$
and $\boldsymbol\omega=dW$ are collinear. This situation is described
by the following lemma, which was used in paper \cite{19}.
\proclaim{Lemma 2.1} If gradient of one smooth function $f(x^1,\ldots,x^n)$
is nonzero in some domain $U\subset\Bbb R^n$ and gradient of another
smooth function $g(x^1,\ldots,x^n)$ is collinear to it in $U$, then
functions $f$ and $g$ are functionally dependent in $U$. This means that
for each point $p\in U$ one can find some neighborhood $O(p)$ and a smooth
function of one variable $\rho(y)$ such that $g=\rho\,\compos f$ in $O(p)$.
\endproclaim
    As an immediate consequence of lemma~2.1 we find that locally in the
neighborhood of each point $q\in\widehat{\Cal M}$ there is some function
$h=h(w)$ such that components of force vector $\bold F$ are determined by
formula \thetag{1.2}. If such function is unique, i\. e\. one for all
points $q$ overall the manifold $\widehat{\Cal M}$, \pagebreak then we say
that {\bf second level of global integration} of structures $\sigma$ and
$s$ is reached. One particular case, when both levels of global integration
are reached, was found in paper \cite{19}. It is described by the following
theorem.
\proclaim{Theorem 7.1} If the section $s$ of factor-bundle $\varOmega
\Cal M=\Cal T\Cal M/U$ corresponding to the force field $\bold F$ of
Newtonian dynamical system admitting the normal shift is nonzero at
all points $q\in\Cal M=M\times\Bbb R^{\sssize +}$, then there is
a global closed 1-form $\boldsymbol\omega$ determining $\bold F$
according to the following formula 
$$
F_k=\frac{N_k}{\omega_{n+1}}-v\sum^n_{i=1}\frac{\omega_i}
{\omega_{n+1}}\,\bigl(2\,N^i\,N_k-\delta^i_k\bigr).\hskip -2em
\tag2.4
$$
\endproclaim
\noindent Formula \thetag{2.4} corresponds to the choice of $h(w)$
being identically equal to unity. Note also that we need not to pass
to universal cover $\widehat{\Cal M}$ in this case. Below we consider
other cases when both levels of integration of structures $\sigma$ and
$s$ are reached, the restriction $s\neq 0$ there is eliminated.
\head
3. Extended tensor fields.
\endhead
     Theorem~1.1 relates formula~2.1 and parameters $b_1,\,\ldots,\,b_n$,
and $a$ in it with structures $\sigma$ and $s$ on Cartesian product
$\Cal M=M\times\Bbb R^{\sssize +}$. However, initially these quantities
were interpreted in a quite different way (see papers \cite{6--16} and
thesis \cite{17}). Function $a$ was interpreted as extended scalar field,
while $b_1,\,\ldots,\,b_n$ are components of extended covectorial field
on $\bold M$. Let's recall appropriate definition.
\definition{Definition 3.1} The function $\bold X$ that to each point
$q=(p,\bold v)$ of tangent bundle $TM$ puts into correspondence some
tensor from the space $T^r_s(p,M)$ at the point $p=\pi(q)$ of $M$ is
called {\bf extended tensor field} of the type $(r,s)$ on $M$.
\enddefinition
    In order to compare note that traditional tensor field $\bold X$
of the type $(r,s)$ on $M$ is a function that maps a point $p$ of $M$,
but not a point of tangent bundle $TM$ as in definition~3.1 above, to
some tensor from the space $T^r_s(p,M)$ at that point $p\in M$. The
idea to extend the concept of tensor field in the sense of
definition~3.1 goes back to Finsler and Cartan (see \cite{20} and
\cite{21}). In book \cite{22} the class of {\bf semibasic tensor
fields}, being subclass in the class of traditional tensor fields
on tangent bundle $TM$, is considered. Theory of {\bf semibasic
tensor fields} constructed in \cite{22} appears to be isomorphic
to the theory of {\bf extended tensor fields} based on definition~3.1.
This fact was discovered by N\.~S\.~Dairbekov when I was making report
in the seminar of Yu\.~G\.~Reshetnyak at the Institute of Mathematics
of Siberian Branch of Russian Academy of Sciences (IM SB RAS) in
October, 2000. Despite to the presence of alternative approach, below
we use extended tensor fields, theory of which was especially developed
(see thesis \cite{17}) for the problems related to Newtonian dynamical
systems in Riemannian and Finslerian manifolds.
\head
4. Norm of covectorial field $\bold b$\\
and global integration of the section $\sigma$.
\endhead
    As we already mentioned above, the quantities $b_1,\,\ldots,\,b_n$
in \thetag{2.1} are interpreted as components of extended covectorial
field $\bold b$ in $\bold M$. But this is extended field of special
form, its components depend on components of velocity vector $\bold v$
in fibers of tangent bundle $TM$ only through their dependence on
modulus of velocity vector $v=|\bold v|$. In paper \cite{23} such fields
were called {\bf fiberwise spherically symmetric}. It is the property of
fiberwise spherical symmetry that gave us the opportunity in \cite{19}
to introduce the manifold $\Cal M=M\times\Bbb R^{\sssize +}$ and find
simple  geometric interpretation of the equations \thetag{1.6} and
\thetag{1.8}. There the field $\bold b$ was associated with the section
$\sigma$ of projectivized cotangent bundle $\Cal P^*\!\Cal M$.\par
     Let's consider the field $\bold b$ in its initial interpretation
and define $|\bold b|$ as a length of covector $\bold b$ in Riemannian
metric of the manifold $M$:
$$
|\bold b|=\sqrt{\shave{\sum^n_{i=1}\sum^n_{j=1}g^{ij}\,b_i\,b_j}}.
$$
The quantity $|\bold b|$ does not depend on local coordinates $x^1,\,
\ldots,\,x^n$, it depends only upon the point $p\in M$ and upon
variable $v\in\Bbb R^{\sssize +}$, which is interpreted as modulus of
velocity vector. Let $f=f(v)$ be a positive function defined on semiaxis
$\Bbb R^{\sssize +}$ and such that the following conditions are fulfilled:
$$
\align
&\quad F(v)=\int\limits^{\,v}_{v_0}\frac{dv}{f(v)}\longrightarrow +\infty
\text{\ \ for \ }v\longrightarrow +\infty,\hskip -2em
\tag4.1\\
&\quad F(v)=\int\limits^{\,v}_{v_0}\frac{dv}{f(v)}\longrightarrow -\infty
\text{\ \ for \ }v\longrightarrow 0.\hskip -2em
\tag4.2
\endalign
$$
Here $v_0$ is some arbitrary positive number from positive semiaxis
$\Bbb R^{\sssize +}$. Let's define $f$-norm of covectorial field
$\bold b$ by the following formula:
$$
\Vert\bold b\Vert=\sup_{p\in M, \ v\in\Bbb R^{\sssize +}}
\left(\frac{|\bold b|}{f(v)}\right).\hskip -2em
\tag4.3
$$
\proclaim{Theorem 4.1} If manifold $M$ is connected and if $f$-norm
of covectorial field $\bold b$ defined by formula \thetag{4.3} is
finite, then first levels of integration of structures $\sigma$ and
$s$ related to the force field of dynamical system is reached.
\endproclaim
    Theorem~4.1 means that if $\Vert\bold b\Vert<\infty$, then the
section $\sigma$ is determined by some global 1-form $\boldsymbol\omega
=dW$. This form can be multivalued, but becomes single-valued upon passage
to universal cover $\widehat{\Cal M}=\widehat M\times\Bbb R^{\sssize
+}$ for the manifold $\Cal M=M\times\Bbb R^{\sssize +}$.\par
    Let's start proving theorem~4.1 by studying problem of local existence
of the form $\boldsymbol\omega$. Thereby we partially resume the content
of paper \cite{19}. Let's choose local coordinates $x^1,\,\ldots,\,x^n,
\,v$ in $\Cal M=M\times\Bbb R^{\sssize +}$. Components of the form
$\boldsymbol\omega$ to be found are bound with components of the field
$\bold b$ by relationships \thetag{1.5}. Therefore in order to find
1-form $\boldsymbol\omega$ it is sufficient to choose proper factor
$\varphi$:
$$
\omega_i=\cases -b_i\,\varphi&\text{for \ }i=1,\,\ldots,\,n,\\
\quad\varphi&\text{for \ }i=n+1.\endcases\hskip -2em
\tag4.4
$$
The condition of closedness of 1-form $\boldsymbol\omega$ is written
in form of relationships
$$
\frac{\partial\omega_i}{\partial x^j}-\frac{\partial\omega_j}
{\partial x^i}=0.\hskip -2em
\tag4.5
$$
Here $v=x^{n+1}$. From \thetag{4.4} and \thetag{4.5} for $i\leqslant n$
and $j\leqslant n$ we derive
$$
\frac{\partial b_i}{\partial x^j}\,\varphi+
\frac{\partial\varphi}{\partial x^j}\,b_i=
\frac{\partial b_j}{\partial x^i}\,\varphi+
\frac{\partial\varphi}{\partial x^i}\,b_j.
\hskip -2em
\tag4.6
$$
From the same relationships \thetag{4.4} and \thetag{4.5} for $i\leqslant n$
and $j=n+1$ we derive
$$
\frac{\partial\varphi}{\partial x^i}=-\frac{\partial b_i}{\partial v}\,
\varphi-\frac{\partial\varphi}{\partial v}\,b_i.\hskip -2em
\tag4.7
$$
Now let's substitute the derivatives $\partial\varphi/\partial x^i$ and
$\partial\varphi/\partial x^j$ calculated according to \thetag{4.7}
into the equations \thetag{4.6}. As a result we get the equations without
entries of $\varphi$. They coincide with \thetag{1.6} exactly. Thus, the
equations \thetag{1.6} form necessary condition for local existence of
closed form $\boldsymbol\omega$ with components \thetag{4.4}. As it was
shown in \cite{19}, these equations constitute sufficient condition as
well. We prove this fact by constructing the solution for the equations
\thetag{4.7}. Let's consider an auxiliary system of partial differential
equations 
$$
\frac{\partial V}{\partial x^i}=b_i(x^1,\ldots,x^n,V)
\text{, \ where \ }i=1,\,\ldots,\,n.\hskip -2em
\tag4.8
$$
This is complete system of Pfaff equations with respect to function
$V(x^1,\ldots,x^n)$. It is compatible. The compatibility condition
for the equations \thetag{4.8} coincides with \thetag{1.6} exactly.
Let's fix some point $p_0\in M$. Without loss of generality we can
assume that local coordinates of the point $p_0$ are equal to zero.
For compatible system of Pfaff equations \thetag{4.8} we set up
the Cauchy problem
$$
V\,\hbox{\vrule height 8pt depth 8pt width 0.5pt}_{\,x^1=\,\ldots\,
\,=\,x^n\,=\,0}=w,\hskip -2em
\tag4.9
$$
where $w>0$. Solution of Cauchy problem \thetag{4.9} exists and is
unique in some neighborhood of the point $p_0$. It is smooth function
of coordinates $x^1,\,\ldots,\,x^n$ and parameter $w$ from right hand
side of \thetag{4.9}:
$$
v=V(x^1,\ldots,x^n,w).\hskip -3em
\tag4.10
$$
For $x^1=\ldots=x^n=0$ due to \thetag{4.9} we get $V(0,\ldots,0,w)=w$.
Therefore
$$
\frac{\partial V}{\partial w}\,\hbox{\vrule height 14pt depth
10pt width 0.5pt}_{\,x^1=\,\ldots\,\,=\,x^n\,=\,0}=1.\hskip -3em
\tag4.11
$$
Let's consider the set of points $q=(p_0,v)$ in $\Cal M$. They constitute
linear ruling in Cartesian product $\Cal M=M\times\Bbb R^{\sssize +}$.
Denote it $l_0=l(p_0)$. The equality \thetag{4.10} means that for any
point $q_0\in l_0$ there is some neighborhood where we have local
coordinates $y^1,\,\ldots,\,y^n,\,w$ related to $x^1,\,\ldots,\,x^n,\,v$
as follows:
$$
\cases x^i=y^i\text{\ \ for \ }i=1,\,\ldots,\,n,\\
v=V(y^1,\ldots,y^n,w).\endcases\hskip -3em
\tag4.12
$$
Inverse passage to $x^1,\,\ldots,\,x^n,\,v$ is determined by the
function $W(x^1,\ldots,x^n,v)$:
$$
\cases y^i=x^i\text{\ \ for \ }i=1,\,\ldots,\,n,\\
w=W(x^1,\ldots,x^n,v).\endcases\hskip -3em
\tag4.13
$$
Function $W(x^1,\ldots,x^n,v)$ is determined in implicit form from
\thetag{4.10} if one treat this equality as an equation with respect
to $w$.\par
    Let's use \thetag{4.12} and \thetag{4.13} in order to simplify the
equations \thetag{4.7}. Instead of function $\varphi(x^1,\ldots,x^n,v)$
in these equations we introduce another function 
$$
\psi(y^1,\ldots,y^n,w)=\varphi(y^1,\ldots,y^n,V(y^1,\ldots,y^n,w)).
\hskip -3em
\tag4.14
$$
The equations \thetag{4.7} are reduced to the following ones with respect
to $\psi$:
$$
\frac{\partial\psi}{\partial y^i}=-B_i\,\psi.\hskip -3em
\tag4.15
$$
The quantities $B_i$ are expressed through partial derivatives of
$V=V(y^1,\ldots,y^n,w)$:
$$
B_i=\frac{1}{Z}\,\frac{\partial Z}{\partial y^i}\text{, \ where \ }
Z=\frac{\partial V}{\partial w}.\hskip -3em
\tag4.16
$$
It's easy to see that \thetag{4.15} is a system of Pfaff equations,
being compatible due to \thetag{4.16}. Moreover, it is explicitly
integrable. General solution of the system of differential equations
\thetag{4.15} has the following form:
$$
\psi=\frac{C(w)}{Z(y^1,\ldots,y^n,w)}.\hskip -3em
\tag4.17
$$
Here $C(w)$ is an arbitrary function of one variable. Now let's use
local invertibility of the relationship \thetag{4.14}:
$$
\varphi(x^1,\ldots,x^n)=\psi(x^1,\ldots,x^n,W(x^1,\ldots,x^n,v)).
\hskip -3em
\tag4.18
$$
From \thetag{4.17} and \thetag{4.18} we can derive general solution
for the equations \thetag{4.7}:
$$
\varphi=C(W)\cdot W_v\text{, \ where \ }W_v=\frac{\partial W}{\partial v}.
\hskip -3em
\tag4.19
$$
Let's turn back to the equations \thetag{4.8} and let's write them with
more details:
$$
b_i(x^1,\ldots,x^1,V(x^1,\ldots,x^n,w))=
\frac{\partial V(x^1,\ldots,x^n,w)}{\partial x^i}.
\hskip -3em
\tag4.20
$$
The relationships \thetag{4.20} are the identities, which are fulfilled
since the function $V(x^1,\ldots,x^n,w)$ is a solution for the system
of equations \thetag{4.8}. Let's substitute $w=W(x^1,\ldots,x^n,v)$ into
\thetag{4.20} in order to express the variable $w$ trough $v$:
$$
b_i(x^1,\ldots,x^1,v)=\frac{\partial V(x^1,\ldots,x^n,w)}
{\partial x^i}\,\hbox{\vrule height 12pt depth 8pt width
0.5pt}_{\,w=W(x^1,\,\ldots,x^n,v)}.\hskip -3em
\tag4.21
$$
If we take into account that $V$ and $W$ determines mutually inverse
changes of variables \thetag{4.12} and \thetag{4.13}, then we can express
right hand side of \thetag{4.21} through partial derivatives of the function
$W(x^1,\ldots,x^n,v)$. This yields
$$
b_i(x^1,\ldots,x^1,v)=-\frac{\nabla_iW}{W_v}\text{, \ where \ }
\nabla_iW=\frac{\partial W}{\partial x^i}.\hskip -3em
\tag4.22
$$
Now let's substitute \thetag{4.22} and \thetag{4.19} into \thetag{4.4}
and calculate components of $\boldsymbol\omega$:
$$
\omega_i=\cases
C(W)\cdot\nabla_iW&\text{for \ }i=1,\,\ldots,\,n,\\
C(W)\cdot W_v&\text{for \ }i=n+1.\endcases\hskip -3em
\tag4.23
$$
If we take the function $C(W)$ being identically equal to unity, then
we get $\boldsymbol\omega=dW$. This means that form $\omega$ just
constructed is closed. However, any other choice of $C(W)$ also yields
closed form $\boldsymbol\omega$.\par
\parshape 17
0pt360pt
0pt 360pt
0pt 360pt
160pt 200pt 160pt 200pt 160pt 200pt 160pt 200pt
160pt 200pt 160pt 200pt 160pt 200pt 160pt 200pt 160pt 200pt
160pt 200pt 160pt 200pt 160pt 200pt 160pt 200pt 
0pt 360pt
     The above method for constructing 1-form $\boldsymbol\omega$ is
purely local yet. The possibility to make it global depends on the
answer to the question --- {\bf how big} is the neighborhood of the point
\vadjust{\vskip 18pt\hbox to 0pt{\kern 8pt\hbox{}\hss}\vskip -18pt}$p_0$,
where the solution of Cauchy problem \thetag{4.9} for the
equations \thetag{4.8} is defined\,? Let $O(p_0)$ be some neighborhood
of the point $p_0$, where such solution does exist, and let $p_1$ be
the point on the boundary of this neighborhood. Let's bind $p_0$ and $p_1$
by a smooth curve $\gamma$ in $M$. It's clear that neighborhood $O(p_0)$
is within the chart where local coordinates $x^1,\,\ldots,\,x^n$ are
defined. Suppose, that the point $p_1$ is also within this chart. Then
curve $\gamma$ can be represented by smooth functions $x^1(t),\ \ldots,
\ x^n(t)$. Coordinates of the point $p_0$ are zero $x^1=\ldots=x^n=0$,
as we have took them above in \thetag{4.9}. Let $t=t_0$ and $t=t_1$
be the values of parameter $t$ for the points $p_0$ and $p_1$ on $\gamma$.
Consider the restriction of the function $V(x^1,\ldots,x^n)$ to $\gamma$:
$$
V(t)=V(x^1(t),\ldots,x^n(t)).\hskip -3em
\tag4.24
$$
Let's differentiate the function \thetag{4.24} with respect to parameter
$t$ and take into account the equations \thetag{4.8}. This yields the
equation
$$
\dot V=\sum^n_{i=1}b_i\,\dot x^i.\hskip -3em
\tag4.25
$$
Since norm \thetag{4.3} is finite, for the derivative $\dot V$ in
\thetag{4.25} we have the estimate:
$$
\left|\frac{\dot V}{f(V)}\right|=\left|\shave{\sum^n_{i=1}}
\frac{b_i\,\dot x^i}{f(V)}\right|\leqslant\frac{|\bold b|}{f(V)}
\cdot|\bold K|\leqslant\Vert\bold b\Vert\cdot|\bold K|.
\hskip -3em
\tag4.26
$$
Here $\bold K$ is a vector with components $\dot x^1,\ldots,\dot x^n$.
It is tangent to $\gamma$. We do not denote it by $\bold v$, since
parameter $t$ on the curve $\gamma$ is not a time. It does not relate
to Newtonian dynamics in \thetag{1.1}. Let $V(t_0)$ and $V(t)$ be
the values of the function \thetag{4.24} at two points on the curve
$\gamma$. Then we have
$$
F(V(t))-F(V(t_0))=\int\limits^{\,t}_{t_0}\frac{\dot V}{f(V)}\,dt.
$$
The function $V(t)$ for $t=t_1$ is not defined. But one can consider
the limit
$$
\lim_{t\to t_1}F(V(t))=F(V(t_0))+\int\limits^{\,\,t_1}_{t_0}
\frac{\dot V}{f(V)}\,dt.\hskip -3em
\tag4.27
$$
Integral in right hand side of \thetag{4.27} is understood as an
improper integral. From the inequality \thetag{4.26} it follows
that such integral absolutely converges:
$$
\int\limits^{\,\,t_1}_{t_0}\left|\frac{\dot V}{f(V)}\right|\,dt
\leqslant\Vert\bold b\Vert\cdot\!\int\limits^{\,\,t_1}_{t_0}|\bold K|
\,dt=\Vert\bold b\Vert\cdot L_\gamma(t_1,t_0).
$$
Here $L_\gamma(t_1,t_0)$ is the length of the segment of curve
$\gamma$ with ending points $p_1$ and $p_0$. It is finite. Hence
the integral in \thetag{4.27} converges, this implies the existence
of finite limit in left hand side of \thetag{4.27}.\par
\parshape 18 0pt 360pt 0pt 360pt 0pt 360pt 160pt 200pt 160pt 200pt
160pt 200pt 160pt 200pt 160pt 200pt 160pt 200pt 160pt 200pt 160pt 200pt
160pt 200pt 160pt 200pt 160pt 200pt 160pt 200pt 160pt 200pt 160pt 200pt
0pt 360pt
    Further let's apply some properties of the function $F(v)$. Let's
remember that it is defined in $\Bbb R^{\sssize +}$, it's monotonic and
increasing, and it satisfies the conditions \thetag{4.1} and \thetag{4.2}.
\vadjust{\vskip 18pt\hbox to 0pt{\kern 8pt\hbox{}\hss}\vskip -18pt}Graph of such
function is drawn on Fig\.~4.2. Due to the above properties of the function
$F(v)$ we can assert that the existence of finite limit in \thetag{4.27}
implies the existence and finiteness of limit
$$
\lim_{t\to t_1}V(t)=\tilde v,\hskip -3em
\tag4.28
$$
the value $\tilde v$ of this limit being positive number from real
semiaxis $\Bbb R^{\sssize +}$:
$$
0<\tilde v<+\infty.
$$
Existence and finiteness of the limit \thetag{4.28} is very important
fact. Now remember that the point $p_1$ corresponding to the value
$t=t_1$ of parameter $t$ is within the chart where local coordinates
$x^1,\,\ldots,\,x^n$ are defined. At the point $p_1$ we can set up the
Cauchy problem similar to \thetag{4.9}:
$$
V\,\hbox{\vrule height 8pt depth 8pt width 0.5pt}_{\,x^1=x^1(p_1),\,
\ldots,\,\,x^n=x^n(p_1)}=\tilde v.\hskip -3em
\tag4.29
$$
The solution of Cauchy problem \thetag{4.29} for Pfaff equations
\thetag{4.8} does exist ant it is unique in some neighborhood $O(p_1)$
of the point $p_1$. It is smooth function of parameter $\tilde v$
and coordinates $x^1,\,\ldots,\,x^n$. Let's denote it $V=\tilde
V(x^1,\ldots,x^n,\tilde v)$. This function also can be restricted to
the curve $\gamma$, where we have the equality
$$
\lim_{t\to t_1}\tilde V(t)=\lim_{t\to t_1}V(t)=\tilde v.\hskip -3em
\tag4.30
$$
From \thetag{4.30} it follows that $V(x^1,\ldots,x^n,w)$ is equal
to $V(x^1,\ldots,x^n,\tilde v)$ on the curve in the region of overlapping
of neighborhoods $O(p_0)$ and $O(p_1)$.
\proclaim{Lemma 4.1} If $f$-norm of covectorial field $\bold b$ is finite,
then the solution of Cauchy problem \thetag{4.9} for the equations
\thetag{4.8} can be continued to any point of chart where local coordinates
$x^1,\,\ldots,\,x^n$ are defined.
\endproclaim
    The equations \thetag{4.8} possess the property of {\bf coordinate
covariance}. This means that their shape doesn't change under the transition
from one set of local coordinates to another. Therefore the solution of
Cauchy problem \thetag{4.9} is a scalar field $V=V(p,w)$ depending on
auxiliary parameter $w$. If $\Vert\bold b\Vert<\infty$ and if we have two
overlapping charts, then scalar field $V$ can be continued from one chart
to another along any curve passing through the region of overlapping.
Now we can strengthen lemma~4.1 as follows.
\proclaim{Lemma 4.2} If $f$-norm of covectorial field $\bold b$ is finite,
then the solution of Cauchy problem \thetag{4.9} for the equations
\thetag{4.8} can be continued to any point $p$ of manifold $M$ along
any curve binding $p$ with the point $p_0$.
\endproclaim
    Note that the result of continuation of scalar field $V$ along the
curve $\gamma$ from the neighborhood of $p_0$ to the point $p$ doesn't
change under continuous deformations of the curve $\gamma$. Therefore
each Cauchy problem for the equations \thetag{4.8} determines some
global scalar field $V$ on universal cover $\widehat M$. There is simple
invariant (non-coordinate) interpretation of scalar field $V$. Indeed,
$V(p,w)$ is a numeric function of the point $p\in\widehat M$ and positive
numeric parameter $w\in\Bbb R^{\sssize +}$, values of this function
also being positive numbers. Its graph is a hypersurface in Cartesian
product $\widehat{\Cal M}=\widehat M\times\Bbb R^{\sssize +}$. It appears
that this hypersurface coincides with integral manifold for involutive
distribution $U=\Ker\sigma$. If we take into account that the point $p_0$,
where Cauchy problem \thetag{4.9} is set up, can be taken for an arbitrary
point in $M$ or, which is more convenient, for an arbitrary point in
universal cover $\widehat M$, then we can reformulate lemma~4.2 as follows.
\proclaim{Lemma 4.3} If $f$-norm of covectorial field $\bold b$ is finite,
then Cartesian product $\widehat{\Cal M}=\widehat M\times\Bbb R^{\sssize +}$
foliates into the disjoint union of integral manifolds of involutive
distribution $U=\Ker\sigma$, each of which being graph for some real-valued
function on $\widehat M$ with the values in $\Bbb R^{\sssize +}$.
\endproclaim
   Thus, the function $v=V(p,w)$ is defined and is single-valued function
on universal cover $\widehat M$ for $M$. Let's define a function $W(p,v)$
such that the conditions $W(p,V(p,w))=w$ and $V(p,W(p,v))=v$ would be
fulfilled.
These conditions are equivalent to requirement that in local coordinates
the changes of variables \thetag{4.12} and \thetag{4.13} are inverse to
each other. Does such function $W(p,v)$ exist on the whole manifold
$\widehat M$\,? The answer to this question depends on solvability of the
equation
$$
v=V(p,w)\hskip -3em
\tag4.31
$$
with respect to variable $w$ at each fixed point $p=p_1$ in $\widehat M$.
In our case, when $\Vert\bold b\Vert<\infty$, the equation \thetag{4.31}
appears to be solvable. Let's prove this fact using lemma~4.3. Suppose
that $v_1\in\Bbb R^{\sssize +}$. Consider the point $q_1=(p_1,v_1)$ of
Cartesian product $\widehat{\Cal M}=\widehat M\times\Bbb R^{\sssize +}$.
Some integral manifold $I$ of distribution $U=\Ker\sigma$ passes through
this point, it is a graph for some function $\psi(p)$. Then $\psi(p_1)
=v_1$. Denote by $w$ the value of this function at the point $p_0$, where
Cauchy problem \thetag{4.9} is set up. This means that we take $w
=\psi(p_0)$. Then submanifold $I$ is a graph for the function $V(p,w)$,
where parameter $w$ is fixed to be equal to $\psi(p_0)$. Therefore
$$
\psi(p)=V(p,\psi(p_0)).\hskip -3em
\tag4.32
$$
Substituting $p=p_1$ into the equality \thetag{4.32}, we get $v_1=V(p_1,
\psi(p_0))$. This means that $w=\psi(p_0)$ is a solution for the equation
$v_1=V(p,w)$ at the point $p=p_1$. Solvability of the equation \thetag{4.31}
means that required function $W(p,v)$ does exist. Now, similar to $V(p,w)$,
it is global, since it is defined on the whole manifold $\widehat M$.\par
    Let's prove that $W(p,v)$ is smooth function. According to the theory
of implicit functions (see \cite{24} or \cite{25}), it is sufficient to
show that the derivative
$$
V_w(p,w)=\frac{\partial V(p,w)}{\partial w}\hskip -3em
\tag4.33
$$
does not vanish. At the point $p=p_0$ the derivative \thetag{4.33} is equal
to unity:
$$
V_w\,\hbox{\vrule height 8pt depth
8pt width 0.5pt}_{\,p=p_0}=1\hskip -3em
\tag4.34
$$
(see relationship \thetag{4.11}). From \thetag{4.8} one can easily derive
the differential equations for the function $V_w$ in local coordinates.
They are the following ones:
$$
\frac{\partial V_w}{\partial x^i}=\frac{\partial b_i}{\partial v}\,V_w.
\hskip -3em
\tag4.35
$$
Similar to $V(p,w)$, the function $V_w$ can be restricted to the curve
$\gamma$. Here we get
$$
W_w(t)=W_w(x^1(t),\ldots,x^n(t),w).\hskip -3em
\tag4.36
$$
For the function \thetag{4.36} from the equations
\thetag{4.35} we derive the differential equation
$$
\frac{dV_w}{dt}=\left(\,\shave{\sum^n_{i=1}}\frac{\partial b_i}
{\partial v}\,\dot x^i\right)\cdot V_w.\hskip -3em
\tag4.37
$$
If the curve $\gamma$ passes through the point $p_0$, then the condition
\thetag{4.34} sets up the Cauchy problem for linear ordinary differential
equation \thetag{4.37}. Its solution does exist and is unique. It is given
by the following formula:
$$
V_w=\exp\left(\,\shave{\int\limits^{\,t}_{t_0}\sum^n_{i=1}}
\frac{\partial b_i}{\partial v}\,\dot x^i\,dt\right)=
\exp\left(\,\shave{\int\limits^{\,t}_{t_0}\bold b_v(\bold K)
\,dt}\right).\hskip -3em
\tag4.38
$$
Here, as in formula \thetag{4.26}, $\bold K$ is the tangent vector
of the curve $\gamma$, its components are $\dot x^1,\ldots,\dot x^n$.
The integral in argument of exponential function in \thetag{4.38} is
a smooth function of parameter of $t$, it has no singular points.
Therefore the value of exponent \thetag{4.38} is nonzero. Hence
$V_w\neq 0$. This provides smoothness of the above function $W(p,v)$.
\pagebreak 
The differential of this function $dW$ is a required 1-form $\boldsymbol
\omega$ (see formula \thetag{4.23} and calculations preceding it).
Thus, theorem~4.1 is proved. This means that under the assumption that
$\Vert\bold b\Vert<\infty$ we reached {\bf first level of global
integration} of geometric structures $\sigma$ and $s$ determining force
field $\bold F$ of Newtonian dynamical system that we consider.
\head
5. Global integration of the section $\sigma$.
\endhead
    Suppose that the condition $\Vert\bold b\Vert<\infty$ is fulfilled.
Let's consider the section $s$ of one-dimensional factor-bundle
$\varOmega\Cal M=\Cal T\Cal M/U$, where $U=\Ker\sigma$. Passing to
universal cover $\widehat{\Cal M}=\widehat M\times\Bbb R^{\sssize +}$
we can consider the section of factor-bundle $\varOmega\widehat{\Cal M}
=\Cal T\widehat{\Cal M}/U$. Due to the condition $\Ker\sigma\nparallel
\Bbb R^{\sssize +}$ (see theorem~1.1 above) such section can be defined
by vectorial field 
$$
\bold X=a(p,v)\cdot\frac{\partial}{\partial v}
$$
in $\widehat{\Cal M}$, or by scalar field $a(p,v)$ in $\widehat{\Cal M}$.
Let's consider the product $\widetilde W=a\cdot W_v$, where $W=W(p,v)$ is
the function which was constructed above in proving theorem~4.1.
Covectors $dW$ and $d\widetilde W$ are collinear and $\boldsymbol\omega=
dW\neq 0$ (see relationship \thetag{2.3} and calculations preceding
lemma~2.1). Let $q_0=(p_0,v_0)$ and $q_1=(p_1,v_1)$ be two points of the
manifold $\widehat{\Cal M}$ lying on the same level hypersurface of the
function $W=W(p,v)$, i\.~e\. such that $W(p_0,v_0)=W(p_1,v_1)$. Suppose
that these points are connected by a curve $\gamma$ lying on the same
level hypersurface as $q_0$ and $q_1$. Then for the difference of
$\widetilde W(p_1,v_1)$ and $\widetilde W(p_0,v_0)$ we get the expression
$$
\widetilde W(p_1,v_1)-\widetilde W(p_0,v_0)=\int\limits^{\,\,t_1}_{t_0}
d\widetilde W(\bold K)\,dt=0.\hskip -2em
\tag5.1
$$
Here $K$ is the tangent vector of curve $\gamma$, while $t_1$ and $t_0$
are the values of parameter $t$ on this curve corresponding to the points
$p_1$ and $p_0$ respectively. Vector $\bold K$ belongs to the kernel of
the form $\boldsymbol\omega=dW$, therefore $dW(\bold K)=0$. Due to
collinearity of covectors $d\widetilde W$ and $dW$ it follows that the
expression $d\widetilde W(\bold K)$ and the integral \thetag{5.1} in whole
do vanish.
\proclaim{Lemma 5.1} If level hypersurfaces of the function $W(p,v)$ are
connected, then $W(p_0,v_0)=W(p_1,v_1)$ implies $\widetilde W(p_0,v_0)=
\widetilde W(p_1,v_1)$.
\endproclaim
     Level hypersurfaces of the function $W(p,v)$ are exactly the integral
manifolds of involutive distribution $U=\Ker\sigma=\Ker\boldsymbol\omega$,
since $\boldsymbol\omega=dW$. Due to lemma~4.3 each of these hypersurfaces
is diffeomorphic to the manifold $\widehat M$. If $M$ is connected, then
$\widehat M$ is also connected. In such situation let's consider the point
$p_0$ where Cauchy problem \thetag{4.9} for the equations \thetag{4.8}
is set up. Here $V(p_0,w)=w$, hence $W(p_0,v)=v$. Let's define the following
function of one variable:
$$
h(v)=\widetilde W(p_0,v).\hskip -2em
\tag5.2
$$
Since $W(p_0,v)=v$, the equality \thetag{5.2} can be rewritten as follows:
$$
\widetilde W(p_0,v)=h(W(p_0,v)).\hskip -2em
\tag5.3
$$
Relying upon lemma~5.1, we can replace $p_0$ in \thetag{5.3} by an
arbitrary point $p$ of universal cover $\widehat M$. Then the equality
\thetag{5.3} looks like
$$
\widetilde W(p,v)=h(W(p,v)).\hskip -2em
\tag5.4
$$
From the equality \thetag{5.4} for the extended scalar field $a$
in \thetag{2.1} we get
$$
a=\frac{h(W)}{W_v}.\hskip -2em
\tag5.5
$$
While components of covectorial field $\bold b$, as we have found above,
are expressed by formula \thetag{4.22}. Substituting \thetag{5.5} and
\thetag{4.22} into the formula \thetag{2.1}, we bring it to the form
\thetag{1.2}. Thereby the functions $h$ and $W$ are now globally
defined for all points $p\in\widehat M$. This means that we reached
{\bf second level of global integration} of geometric structures $\sigma$
and $s$.
\proclaim{Theorem 5.1} If Riemannian manifold $M$ is connected and if
$f$-norm of covectorial field $\bold b$ is finite, then for geometric
structures determining force field of Newtonian dynamical system admitting
the normal shift of hypersurfaces in $M$ both levels of global integration
are reached.
\endproclaim
\head
6. Monodromy transformations.
\endhead
    Force field $F$ and geometric structures $\sigma$ and $s$ determining
this field are related to the manifold $\Cal M$. However, in integrating
these structures we are to pass to universal cover $\widehat M$. Therefore
the functions $W$ and $h$, which were constructed above, should contain
a discrete symmetry determined by first fundamental group $\pi_1(M)$.
Group $\pi_1(M)$ acts in $\widehat M$ by discrete transformations, and
$M$ coincides with the result of factorization of $\widehat M$ with respect
to such action: $M=\widehat M/\pi_1(M)$. Take $g\in\pi_1(M)$. Let's compare
two functions $V(p,w)$ and $V(g^{-1}(p),v)$. If we localize them in the
neighborhood of the point $p_0$, where Cauchy problem \thetag{4.9} is set
up, then, upon passing from $\widehat M$ to $M$ by canonical projection,
these functions appears to be the solutions of the same system of
differential equations \thetag{4.8}. Now let's consider the following
function of one variable:
$$
\rho(w)=V(g^{-1}(p_0),w).
$$
Then for the values of functions $V(p,w)$ and $V(g^{-1}(p),v)$ at the
point $p_0$ we get
$$
\xalignat 2
&V(p,w)\,\hbox{\vrule height 8pt depth 8pt width 0.5pt}_{\,p=p_0}=w,
&&V(g^{-1}(p),w)\,\hbox{\vrule height 8pt depth 8pt width 0.5pt}_{\,p=p_0}=
\rho(w).\hskip -2em
\tag6.1
\endxalignat
$$
Each of the relationships \thetag{6.1} can be treated as Cauchy problem
for the equations \thetag{4.8} in the neighborhood of the point $p_0$.
Due to the uniqueness of solutions of such Cauchy problems, we get the
following relationship:
$$
V(g^{-1}(p),w)=V(p,\rho(w)).\hskip -2em
\tag6.2
$$
Initially the relationship \thetag{6.2} is fulfilled in some neighborhood
of the point $p_0$ on $M$. However, since we can continue functions
$V(p,w)$ and $V(g^{-1}(p),v)$ along any curve $\gamma$ in $M$, it is
fulfilled as an identity on universal cover $\widehat M$. Function $W(p,v)$
is defined as the solution of the equation $v=V(p,w)$ with respect to $w$
for fixed $p$ (see above). Therefore from \thetag{6.2} we derive
$$
\rho(W(g^{-1}(p),v))=W(p,v).\hskip -2em
\tag6.3
$$
The relationship \thetag{6.3} can be rewritten as follows:
$$
W(g(p),v)=\rho(W(p,v)).\hskip -2em
\tag6.4
$$
So, each element $g$ from first fundamental group $\pi_1(M)$ appears to
be related to some function $\rho=\rho_g(w)$. From \thetag{6.4} it's
easy to derive the relationship 
$$
\rho_{g_1\cdot\,g_2}=\rho_{g_1}\compos\rho_{g_1}.
$$
This means that we have a representation of the group $\pi_1(M)$ by
transformations of real semiaxis $\Bbb R^{\sssize +}$ given by smooth
strictly monotonic increasing functions $\rho_g(w)$. Such
transformations are usually called {\bf monodromy transformations}.
\par
    Extended scalar field $a(p,v)$ on universal cover $\widehat M$
is obtained by lifting the corresponding scalar field $a$ from $M$.
Therefore
$$
a(g(p),v)=a(p,v).\hskip -2em
\tag6.5
$$
Let's differentiate the equality \thetag{6.4} with respect to $v$ for
fixed $p\in M$. This yields 
$$
W_v(g(p),v)=\rho'(W(p,v))\cdot W_v(p,v).\hskip -2em
\tag6.6
$$
If we take into account \thetag{6.5} and \thetag{6.6}, then formula
\thetag{5.5} can be written as
$$
a=\frac{h(\rho^{-1}(W^\bullet))\cdot\rho'(\rho^{-1}(W^\bullet))}
{W^\bullet_v}\text{, \ where \ }W^\bullet=W(g(p),v).\hskip -2em
\tag6.7
$$
Due to the relationships \thetag{6.4} and \thetag{6.7} we can associate
each element $g$ of first fundamental group $\pi_1(M)$ with transformations
of the form 
$$
\aligned
&W(p,v)\longrightarrow \rho(W(p,v)),\\
&h(w)\longrightarrow h(\rho^{-1}(w))\cdot\rho'(\rho^{-1}(w)),
\endaligned\hskip -2em
\tag6.8
$$
where $\rho=\rho_g$. In the framework of local approach transformations
of the form \thetag{6.8} were obtained in thesis \cite{17} (see \S\,5
in Chapter~\uppercase\expandafter{\romannumeral 7}) as transformations
changing the pair of functions $(h,W)$, but not changing the force field
$\bold F$ given by formula \thetag{1.2}. They were called {\bf gauge
transformations}.
\proclaim{Theorem 6.1} Suppose that $M$ is connected Riemannian manifold
equipped with Newtonian dynamical system admitting the normal shift of
hypersurfaces. In this situation if $f$-norm of covectorial field $\bold
b$ corresponding to force field $\bold F$ of such system is finite, then
\roster
\rosteritemwd=0pt
\item"1)" the set of pairs of functions $(h,W)$ determining
$\bold F$ by formula \thetag{1.2} globally at all points of
$\widehat M$ is not \pagebreak empty;
\item"2)" group $\pi_1(M)$ acts in the set of such pairs of functions
by means of gauge transformations of the form \thetag{6.8}.
\endroster
\endproclaim
    First proposition in theorem~6.1 is direct consequence of theorem~5.1.
Second proposition of this theorem was proved above.
\head
7. Second problem of globalization.
\endhead
     Let $S$ be a hypersurface in $M$, and suppose that $p\in S$. Consider
the following initial data for the system of equations \thetag{1.1}:
$$
\xalignat 2
&\quad x^k\,\hbox{\vrule height 8pt depth 8pt width 0.5pt}_{\,t=0}
=x^k(p),
&&\dot x^k\,\hbox{\vrule height 8pt depth 8pt width 0.5pt}_{\,t=0}
=\nu(p)\cdot n^k(p).\hskip -2em
\tag7.1
\endxalignat
$$
Here $n^k(p)$ are components of unitary normal vector $\bold n$ to
$S$ at the point $p$. Initial data \thetag{7.1} define the trajectory
of dynamical system \thetag{1.1} coming out from the point $p$ in
the direction of normal vector $\bold n(p)$, while the quantity
$\nu(p)$ in \thetag{7.1} determines modulus of initial velocity for
such trajectory.\par
    Let's choose and fix some point $p_0\in S$, then consider a smooth
function $\nu(p)$ defined on $S$ in some neighborhood of the point $p_0$.
Suppose that 
$$
\nu(p_0)=\nu_0>0.\hskip -2em
\tag7.2
$$
Then in some (possibly smaller) neighborhood of the point $p_0$ the
function $\nu(p)$ is positive. Restricting $\nu(p)$ to such neighborhood,
we use it for to determine initial velocity in \thetag{7.1}. As a result
we get the whole family of trajectories of dynamical system \thetag{1.1}.
The displacement of points of hypersurface $S$ along such trajectories
determines shift maps $f_t\!:S'\to S'_t$. Relying upon the theorem on
existence, uniqueness, and smooth dependence on initial data for the
systems of ordinary differential equations (see \cite{26} and \cite{27}),
we can take shift maps $f_t\!:S'\to S'_t$ to be defined in some neighborhood
$S'$ of the point $p_0$ on $S$ for all values of parameter $t$ from some
interval $(-\varepsilon,\,+\varepsilon)$ on real axis. At the expense of
further restriction of neighborhood $S'$ and the interval $(-\varepsilon,
\,+\varepsilon)$ one can achieve the situation, when shift maps would
become diffeomorphisms, while their images $S'_t$ would become smooth
hypersurfaces, disjoint union of which would fill some neighborhood
of the point $p_0$ in $M$. Moreover, at the expense of restricting
the neighborhood $S'$ and the interval $(-\varepsilon,\,+\varepsilon)$
one can achieve the transversality of hypersurfaces $S_t$ and shift
trajectories at all pints of their intersection.
\definition{Definition 7.1} Shift $f_t\!:S'\to S'_t$ of a part $S'$ of
hypersurface $S$ along trajectories of Newtonian dynamical
system \thetag{1.1} is called {\bf a normal shift} if all hypersurfaces
$S'_t$ arising in the process of shifting are perpendicular to
shift trajectories.
\enddefinition
\definition{Definition 7.2} Newtonian dynamical system \thetag{1.1}
with force field $\bold F$ is called a system {\bf admitting normal
shift} in strong sense\footnote{First we used the definition without
normalizing condition \thetag{7.2} for the function $\nu(p)$. Such
definition was called the normality condition. Definition~7.2
strengthens this condition making it more restrictive with respect
to the choice of force field $\bold F$ of dynamical system \thetag{1.1}.
Therefore it is called strong normality condition.}\ if for any
hypersurface $S$ in $M$, for any point $p_0\in S$, \pagebreak and for
any real number $\nu_0>0$ there exists a neighborhood $S'$ of the point
$p_0$ on $S$ and there exits a smooth positive in $S'$ function $\nu(p)$
normalized by the condition \thetag{7.2}, and such that the shift
$f_t\!:S'\to S'_t$ defined by this function is a normal shift in the
sense of definition~7.1.
\enddefinition
\adjustfootnotemark{-1}
     Definitions~1.1 and 1.2 appeared to be very fruitful. On the base
of these definitions in papers \cite{1--16} the theory of dynamical
systems admitting the normal shift was constructed. However, in these
definitions we observe the series restrictions making theory very local.
The most displeasing is the necessity to replace whole hypersurface by
by a neighborhood $S'$ of marked point $p_0$. So we meet the problem of
finding situations, when one could provide the possibility to define a
function $\nu(p)$ and shift maps $f_t\!:S\to S_t$ globally on the whole
hypersurface $S$. This problem was called {\bf a second problem of
globalization}. It was formulated by A.~S.~Mishchenko when I was reporting
the results of thesis \cite{17} and succeeding papers \cite{23}, \cite{28},
and \cite{29} in his seminar at Moscow State University.\par
    Note that second problem of globalization is closely related to the
{\bf first problem of globalization}, which was considered in paper
\cite{19}. First problem of globalization was formulated by S.~E.~Kozlov
and Yu\.~R\.~Romanovsky when we were discussing the results of thesis
\cite{17} and succeeding papers \cite{23}, \cite{28}, and \cite{29} in
the seminar of N.~Yu\.~Netsvetaev at Saint-Petersburg department of
Steklov Mathematical Institute. 
\head
8. Choosing initial velocity\\
in the construction of normal shift.
\endhead
     Let $M$ be connected Riemannian manifold equipped with a Newtonian
dynamical system \thetag{1.1} admitting the normal shift of hypersurfaces.
Suppose that $f$-norm of covectorial field $\bold b$ corresponding to
the force field $\bold F$ of this system is finite. Let's choose and fix
some hypersurface $S$ and some point $p_0$ on it. We choose local
coordinates $u^1,\,\ldots,\,u^{n-1}$ on $S$ in some neighborhood of the
$p_0$ and local coordinates $x^1,\,\ldots,\,x^n$ in the manifold $M$ in
a neighborhood of the same point $p_0$. Without loss of generality one
can assume that coordinates of the point $p_0$ are zero: $u^1=\ldots
=u^{n-1}=0$ and $x^1=\ldots=x^n=0$. Now hypersurface $S$ in a neighborhood
of the point $p_0$ can be represented parametrically by the following
functions:
$$
\align
x^1&=x^1(u^1,\ldots,u^{n-1}),\hskip -2em\\
.\ .\ &.\ .\ .\ .\ .\ .\ .\ .\ .\ .\ .\ .\ .\ .\ .
\hskip -2em\tag8.1\\
x^n&=x^n(u^1,\ldots,u^{n-1}).\hskip -2em
\endalign
$$
The choice of local coordinates $u^1,\,\ldots,\,u^{n-1}$ determines
coordinate tangent vectors $\boldsymbol\tau_1,\,\ldots,\,\boldsymbol
\tau_{n-1}$ forming a base in tangent hyperplane to hypersurface $S$.
They can be determined by the relationships 
$$
\boldsymbol\tau_k=\sum^n_{i=1}\frac{\partial x^i}{\partial u^k}
\cdot\frac{\partial}{\partial x^i}\text{, \ where \ }k=1,\,\ldots,
\,n-1.\hskip -2em
\tag8.2
$$\par
     Second problem of globalization is related to the problem of
constructing smooth positive function $\nu(p)$ on $S$ which would
be normalized by the condition \thetag{7.2} and would define the
normal shift of hypersurface $S$ along trajectories of dynamical
system \thetag{1.1} by fixing the value of initial velocity in 
\thetag{7.1}. In the framework of local approach an algorithm of
constructing such function $\nu(p)$ was found in \cite{3} and
\cite{7} (see also Chapter~\uppercase\expandafter{\romannumeral 5}
of thesis \cite{17}). Omitting details, we shall only use the fact
that $\nu(p)=\nu(u^1,\ldots,u^n)$ is constructed as a solution of
the equations 
$$
\frac{\partial\nu}{\partial u^k}=-\frac{(\bold F\,|\,\boldsymbol
\tau_k)}{\nu}\text{, \ where \ }k=1,\,\ldots,\,n-1.\hskip -2em
\tag8.3
$$
In calculating scalar product $(\bold F\,|\,\boldsymbol\tau_k)$ in
\thetag{8.3} now we can use explicit formula \thetag{2.1} for
components of force field $\bold F$. Moreover, let's take into
account the relationships \thetag{8.2} which determine components
of vectors $\boldsymbol\tau_1,\,\ldots,\,\boldsymbol\tau_{n-1}$.
This yields
$$
\frac{\partial\nu}{\partial u^k}=\sum^n_{i=1}b_i(x^1,\ldots,x^n,\nu)
\,\frac{\partial x^i}{\partial u^k}\hskip -2em
\tag8.4
$$
where $k=1,\,\ldots,\,n-1$. In deriving \thetag{8.4} we take into account
that at initial instant of time $t=0$ the velocity vector $\bold v$ is
directed along normal vector to $S$ (see initial data \thetag{7.1}).
Therefore $\bold v\perp\boldsymbol\tau_k$ and $\bold N\perp\boldsymbol
\tau_k$. Modulus of velocity vector for $t=0$ coincides with $\nu$. The
dependence of $x^1,\,\ldots,\,x^n$ on $u^1,\,\ldots,\,u^{n-1}$ in the
equations \thetag{8.4} is determined by functions \thetag{8.1}. The same
functions determine partial derivatives $\partial x^i/\partial u^k$ in
right hand side of these equations.\par
     The equations \thetag{8.4} form complete system of Pfaff equations
for the function $\nu$. It is compatible. Its compatibility follows from
the relationships \thetag{1.6}. Normalizing condition \thetag{7.2} sets
up the Cauchy problem for Pfaff equations \thetag{8.4}. Such Cauchy
problem has unique solution in some neighborhood of the point $p_0$ on
$S$. Now we are to study whether it's possible to continue this solution
to whole hypersurface $S$. Let's compare the equations \thetag{8.4} with
the equations \thetag{4.8}, for which the Cauchy problem \thetag{4.9} at
the point $p_0$ is set up. The equations \thetag{8.4} can be treated as
the restrictions of the equations \thetag{4.8} from $M$ to $S$. If $V(x^1,
\ldots,x^n,w)$ is the solution of Cauchy problem \thetag{4.9} for the
equations \thetag{4.8}, then, substituting $w=\nu_0$ and substituting the
functions \thetag{8.1} for $x^1,\,\ldots,\,x^n$, we get the solution of
Cauchy problem \thetag{7.2} for the equations \thetag{8.4}. This fact
indicates the way for solving second problem of globalization.
\proclaim{Theorem 8.1} Suppose that $M$ is connected Riemannian manifold
equipped with Newtonian dynamical system admitting the normal shift of
hypersurfaces. In this situation if $f$-norm of covectorial field $\bold
b$ corresponding to force field $\bold F$ of such system is finite, then
for any hypersurface $S$ in $M$ there is a function $\nu(p)$ normalized
by the condition \thetag{7.2} such that it determines modulus of initial
velocity in the construction of normal shift for $S$. This function is
continued globally to any point $p\in S$ along any curve lying on $S$,
though thereby it may appear to be multivalued.
\endproclaim
     The multivalued function $\nu(p)$ may arise since function $V(p,w)$
is defined not in $M$, but in universal cover $\widehat M$. Remember that
first fundamental group $\pi_1(M)$ acts as a group of discrete
transformations in $\widehat M$. Let's define the following subgroup:
$$
G_{\bold F}=\{g\in\pi_1(M)\text{\ \ such that \ }
W(g(p),w)\equiv W(p,w)\}.
$$
Subgroup $G_{\bold F}$ is a characteristic (topological invariant) of
the force field $\bold F$ in $\bold M$. It is formed by elements
monodromy transformations $\rho_g$ \pagebreak for which are identical.
\par
    Let $S$ be a hypersurface in $M$. It is known that the immersion
$S\subset M$ determines homomorphism of fundamental groups $\pi_1(S)
\to\pi_1(M)$.
\proclaim{Theorem 8.2} Under the assumption of theorem~8.1 the function
$\nu$ on $S$ is single-valued if and only if the image of the group
$\pi_1(S)$ under the immersion homomorphism $\pi_1(S)\to\pi_1(M)$ is
contained in subgroup $G_{\bold F}$.
\endproclaim
    Note that if $M$ is simply connected or hypersurface $S$ is simply
connected, then the condition providing univalence of $\nu$ is fulfilled.
\head
9. Acknowledgements.
\endhead
    I am grateful to A.~S.~Mishchenko for the invitation to visit
Moscow and for the opportunity to report the results of thesis \cite{17}
and succeeding papers \cite{23}, \cite{28}, and \cite{29} in his
seminar at Moscow State University. I am grateful to N.~Yu\.~Netsvetaev
for the invitation to visit Saint-Petersburg and for the opportunity
to report the same results in the seminar at Saint-Petersburg department
of Steklov Mathematical Institute. I am grateful to all participants of
both seminars mentioned above and to my colleague E.~G.~Neufeld from
Bashkir State University for fruitful discussions which stimulated
preparing this paper.\par
     This work is supported by grant from Russian Fund for Basic Research
(project No\nolinebreak\.~00\nolinebreak-01-00068, coordinator 
Ya\.~T.~Sultanaev), and by grant from Academy of Sciences of the
Republic Bashkortostan (coordinator N.~M.~Asadullin). I am grateful
to these organizations for financial support.\par


\newpage
\Refs
\ref\no 1\by Boldin~A\.~Yu\., Sharipov~R\.~A\.\book Dynamical systems
accepting the normal shift\publ Preprint No\.~0001-M of Bashkir State
University\publaddr Ufa\yr April, 1993
\endref
\ref\no 2\by Boldin~A.~Yu\., Sharipov~R.~A.\paper Dynamical systems
accepting the normal shift\jour Theoretical and Mathematical Physics (TMF)
\vol 97\issue 3\yr 1993\pages 386--395\moreref see also chao-dyn/9403003
in Electronic Archive at LANL\footnotemark
\endref
\footnotetext{Electronic Archive at Los Alamos national Laboratory of USA
(LANL). Archive is accessible through Internet 
{\bf http:/\negskp/xxx\.lanl\.gov}, it has mirror site 
{\bf http:/\negskp/xxx\.itep\.ru} at the Institute for Theoretical and
Experimental Physics (ITEP, Moscow).}
\ref\no 3\by Boldin~A.~Yu\., Sharipov~R.~A.\paper Multidimensional
dynamical systems accepting the normal shift\jour Theoretical and
Mathematical Physics (TMF)\vol 100\issue 2\yr 1994\pages 264--269
\moreref see also patt-sol/9404001 in Electronic Archive at LANL
\endref
\ref\no 4\by Boldin~A.~Yu\., Sharipov~R.~A.\paper Dynamical systems
accepting the normal shift\jour Reports of Russian Academy of Sciences
(Dokladi RAN)\vol 334\yr 1994\issue 2\pages 165--167
\endref
\ref\no 5\by Sharipov~R.~A.\paper Problem of metrizability for
the dynamical systems accepting the normal shift\jour Theoretical and
Mathematical Physics (TMF)\yr 1994\vol 101\issue 1\pages 85--93\moreref
see also solv-int/9404003 in Electronic Archive at LANL
\endref
\ref\no 6\by Boldin~A.~Yu\., Dmitrieva~V.~V., Safin~S.~S., Sharipov~R.~A.
\paper Dynamical systems accepting the normal shift on an arbitrary 
Riemannian manifold\jour Theoretical and Mathematical Physics (TMF)
\yr 1995\vol 105\issue 2\pages 256--266\moreref\inbook see also
``{Dynamical systems accepting the normal shift}'', Collection of papers
\publ Bashkir State University\publaddr Ufa\yr 1994\pages 4--19
\moreref see also hep-th/9405021 in Electronic Archive at LANL
\endref
\ref\no 7\by Boldin~A.~Yu\., Bronnikov~A.~A., Dmitrieva~V.~V.,
Sharipov~R.~A.\paper Complete normality conditions for the dynamical
systems on Riemannian manifolds\jour Theoretical and Mathematical
Physics (TMF)\yr 1995\vol 103\issue 2\pages 267--275\moreref\inbook
see also ``{Dynamical systems accepting the normal shift}'', Collection
of papers\publ Bashkir State University\publaddr Ufa\yr 1994
\pages 20--30\moreref see also astro-ph/9405049 in Electronic Archive
at LANL
\endref
\ref\no 8\by Boldin~A\.~Yu\.\paper On the self-similar solutions of 
normality equation in two-dimensional case\inbook ``{Dynamical systems
accepting the normal shift}'', Collection of papers\publ Bashkir State
University\publaddr Ufa\yr 1994\pages 31--39\moreref see also
patt-sol/9407002 in Electronic Archive at LANL
\endref
\ref\no 9\by Sharipov~R.~A.\paper Metrizability by means of conformally
equivalent metric for the dynamical systems\jour Theoretical and
Mathematical Physics (TMF)\yr 1995\vol 105\issue 2\pages 276--282
\moreref\inbook see also ``{Integrability in dynamical systems}''\publ
Institute of Mathematics, Bashkir Scientific Center of Ural branch of
Russian Academy of Sciences (BNC UrO RAN)\publaddr Ufa\yr 1994
\pages 80--90
\endref
\ref\no 10\by Sharipov~R\.~A\.\paper Dynamical systems accepting normal
shift in Finslerian geometry,\yr November, 1993\finalinfo 
unpublished\footnotemark
\endref
\footnotetext{Papers \cite{3--18} are arranged here in the order they
were written. However, the order of publication not always coincides with
the order of writing.}\adjustfootnotemark{-2}
\ref\no 11\by Sharipov~R\.~A\.\paper Normality conditions and affine
variations of connection on Riemannian manifolds,\yr December, 1993
\finalinfo unpublished
\endref
\ref\no 12\by Sharipov~R.~A.\paper Dynamical system accepting the normal
shift (report at the conference)\jour see in Progress in Mathematical
Sciences (Uspehi Mat\. Nauk)\vol 49\yr 1994\issue 4\page 105
\endref
\ref\no 13\by Sharipov~R.~A.\paper Higher dynamical systems accepting 
the normal shift\inbook ``{Dynamical systems accepting the normal 
shift}'', Collection of papers\publ Bashkir State University\publaddr 
Ufa\yr 1994\linebreak\pages 41--65
\endref
\ref\no 14\by Dmitrieva~V.~V.\paper On the equivalence of two forms
of normality equations in $\Bbb R^n$\inbook ``{Integrability in dynamical
systems}''\publ Institute of Mathematics, Bashkir Scientific Center of
Ural branch of Russian Academy of Sciences (BNC UrO RAN)\publaddr
Ufa\yr 1994\pages 5--16
\endref
\ref\no 15\by Bronnikov~A.~A., Sharipov~R.~A.\paper Axially
symmetric dynamical systems accep\-ting the normal shift in $\Bbb R^n$
\inbook ``{Integrability in dynamical systems}''\publ Institute of
Mathematics, Bashkir Scientific Center of Ural branch of Russian Academy
of Sciences (BNC UrO RAN)\publaddr Ufa\yr 1994\linebreak\pages 62--69
\endref
\ref\no 16\by Boldin~A.~Yu\., Sharipov~R.~A.\paper On the solution
of normality equations in the dimension $n\geqslant 3$\jour Algebra and
Analysis (Algebra i Analiz)\vol 10\yr 1998\issue 4\pages 37--62\moreref
see also solv-int/9610006 in Electronic Archive at LANL
\endref
\ref\no 17\by Sharipov~R.~A.\book Dynamical systems admitting the normal
shift\publ Thesis for the degree of Doctor of Sciences in Russia\yr 1999
\moreref English version of thesis is submitted to Electronic Archive at 
LANL, see archive file math.DG/0002202 in the section of Differential 
Geometry\footnotemark
\endref
\footnotetext{For the convenience of reader we give direct reference
to archive file. This is the following URL address:
{\bf http:/\negskp/xxx\.lanl\.gov/eprint/math\.DG/0002202}\,.}
\ref\no 18\by Boldin~A.~Yu\.\book Two-dimensional dynamical systems
admitting the normal shift\publ Thesis for the degree of Candidate of
Sciences in Russia\yr 2000\moreref English version of thesis is
submitted to Electronic Archive at LANL, see archive file math.DG/0011134
in the section of Differential Geometry
\endref
\ref\no 19\by Sharipov~R.~A.\paper First problem of globalization
in the theory of dynamical systems admitting the normal shift of
hypersurfaces\jour Paper math.DG math.DG/0101150 in Electronic
Archive at LANL\yr 2001
\endref
\ref\no 20\by Finsler\book \"Uber Kurven and Flachen in algemeinen Raumen
\publ Dissertation\publaddr G\"ottin\-gen\yr 1918
\endref
\ref\no 21\by Cartan~E\. \book Les espaces de Finsler\publ Actualites 79
\publaddr Paris\yr 1934
\endref
\ref\no 22\by Sharafutdinov~V.~A.\book Integral geometry of tensor
fields\publ VSP\publaddr Utrecht, The Netherlands\yr 1994
\endref
\ref\no 23\by Sharipov~R.~A.\paper Newtonian normal shift in
multidimensional Riemannian geometry\jour Paper math.DG/0006125
in Electronic Archive at LANL\yr 2000
\endref
\ref\no 24\by  Kudryavtsev~L.~D.\book  Course  of  mathematical 
analysis, Vol\.~\uppercase\expandafter{\romannumeral 1},
\uppercase\expandafter{\romannumeral 2}\publ ``Nauka'' publishers
\publaddr Mos\-cow\nolinebreak\yr 1985
\endref
\ref\no 25\by Ilyin~V.~A., Sadovnichiy~V.~A., Sendov~B.~H.
\book Mathematical analysis \publ ``Nauka'' publishers\publaddr 
Moscow\yr 1979
\endref
\ref\no 26\by Petrovsky~I.~G.\book Lectures on the theory of ordinary
differential equations\publ Moscow State University publishers\yr 1984
\publaddr Moscow
\endref
\ref\no 27\by Fedoryuk~M.~V.\book Ordinary differential equations
\yr 1980\publ ``Nauka'' publishers\publaddr Moscow
\endref
\ref\no 28\by Sharipov~R.~A.\paper Newtonian dynamical systems
admitting normal blow-up of points\jour Paper math.DG/0008081
in Electronic Archive at LANL\yr 2000
\endref
\ref\no 29\by Sharipov~R.~A.\paper On the solutions of weak normality
equations in multidimensional case\jour Paper math.DG/0012110
in Electronic Archive at LANL\yr 2000
\endref
\endRefs
\enddocument
\end